\documentclass[a4paper]{amsart}
\usepackage[T1]{fontenc}
\usepackage{hhline}
\usepackage{longtable}
\usepackage{amscd}
\usepackage{epsfig}
\usepackage{amsmath}
\usepackage{psfrag}
\usepackage{graphicx}
\def\R {I\!\!R}
\begin{document}
\title[operational aircraft parameters Reducing Noise Emission]{Optimization of operational aircraft parameters Reducing Noise Emission}
\author []{Lina Abdallah, Mounir
Haddou and Salah Khardi  \\
\\}

\begin{abstract}
The objective of this paper is to develop a model and a minimization
method to provide flight path optimums reducing aircraft noise in
the vicinity of airports. Optimization algorithm has solved a
complex optimal control problem, and generates flight paths
minimizing aircraft noise levels. Operational and safety constraints
have been considered and their limits satisfied. Results are here
presented and discussed.
\end{abstract}
\maketitle {\bf{Nomenclature}}
\begin{displaymath}
\begin{array}{llllll}
\hline\\
C_{z_{\alpha}}&& \textrm{slope of lift coefficient
curve}\hspace{2cm}&\delta_x&&\textrm{throttle setting}\\
C_{x_0}&&\textrm{drag coefficient}\hspace{2cm}&k_i&&\textrm{induced drag parameter}\\
T&&\textrm{thrust}, N\hspace{2cm}&V&&\textrm{speed of aircraft}, m/s\\
L&&\textrm{lift, N}\hspace{2cm}&S&&\textrm{wing area}, m^2\\
D&&\textrm{drag, N}\hspace{2cm}&c&&\textrm{speed of sound}, m/s\\
g&&9.8\,m/s^2\hspace{2cm}&M&&\textrm{mach number}, V/c\\
x&&\textrm{horizontal distance}, m\hspace{2cm}&\rho&&\textrm{air density}, kg/m^3\\
y&&\textrm{lateral distance}, m \hspace{2cm}&d&&\textrm{nozzle diameter}, m\\
h&&\textrm{aircraft height}, m\hspace{2cm}& s&&\textrm{area of
coaxial engine nozzle}, m^2\\
k&&\textrm{induced drag coefficient}\hspace{2cm}& v&&\textrm{speed of gas}, m/s\\
t&&\textrm{time}, s\hspace{2cm}& w&&\textrm{density exponent}\\
\mu&&\textrm{roll angle}, rad \hspace{2cm}& \tau&&\textrm{temperature, K}\\
\alpha&&\textrm{angle of attack}, rad\hspace{2cm}&\theta
&&\textrm{directivity angle}, rad \\
\gamma&&\textrm{flight path angle}, rad\hspace{2cm}&R&&\textrm{source-to-observer distance}, m\\
\chi&&\textrm{yaw angle}, rad\hspace{2cm}& m &&\textrm{aircraft mass}, kg \\
\\
 \hline\\ \\
\end{array}
\end{displaymath}

\section{Introduction}
Since the introduction of jet aircraft in the 1960, aircraft noise
produced in the vicinity of airports has represented a serious
social and environmental issue. It is a continuing source of
annoyance in nearby communities. The importance of that problem has
been highlighted by the increased public concern for environmental\\
issues. To deal with this problem, aircraft manufacturers and public
establishments are engaged in research on technical and theoretical
approaches for noise reduction concepts that should be applied to
new aircraft. The ability to assess noise exposure accurately is an
increasingly important factor in the design and implementation of
any airport improvements.
\\ \\
Aircraft are complex noise sources and the emitted intensities vary
with the type of aircraft, in particular, with the type of engines
and with the implemented flight procedures. The noise contour
assessment due to the variety of flight route schemes and predicted
procedures is also complex. A set of data must be used which
includes noise data, flight path parameters and their features, and
environmental conditions affecting outdoor sound propagation. Three
development initiatives are available for the reduction of aircraft
noise: (1) innovative passive technologies required by the industry
for developing environmentally compatible and economically viable
aircraft, (2) advanced active technologies such as computational
aeroacoustics, active control, advanced propagation and prediction
methods, (3) reliable trajectory and procedures optimization which
can be used to determine optimal landing approach for any arbitrary
aircraft at any given airport. The last action will be particularly
emphasized in the next sections.
\\ \\
A number of calculation programs of aircraft noise impact have been
developed over the last 30 years. They have been widely used by
aircraft manufacturers and airport authorities.
Their reliability and results efficiency to assess the
real impact of aircraft noise have not been proved conclusively.
They are complex, very slow and can not be planned for on-line and
on-site use. That is the reason why, the model described in this
paper, generating optimal trajectories minimizing noise, is
considered as a promising
scientific plan.\\
\\
Several optimization codes for NLP exist in the literature. After a
large number of test and comparisons, we choose KNITRO \cite{Kni}
which known for its performances and robustness, this software is
efficient to solve general nonlinear programming. We will explain in
the following sections how the considered optimal control problem is
discretized and solved. Numerical results have been analyzed and
their reliability and flexibility have been proved. We have
demonstrated the effectiveness of computation and its application to
aircraft noise reduction. The objective of that alternative research
is to develop high payoff models to enable a safe, and
environmentally compatible and economical aircraft. We should make
large profits in terms of noise abatement in comparison with the
expected noise control systems in progress. These systems, which are
not in an advanced step, in particular at low frequencies, are still
ineffective or impractical. Actually, the low-frequency broadband
generated by the engines represents a significant source of
environmental noise. Their radiation during flight operations is
extremely difficult to attenuate using the mentioned systems and is
capable of propagating over long distance \cite{Smith96}.\\
\\
Details of trajectory and aircraft noise models, and optimal control
problems are presented in section 2 and 3 while the last section is
devoted to numerical experiments.
\section{Optimal Control Problem}
\subsection{Equation of Motion}
In general, the system of differential equations commonly employed
in aircraft trajectory analysis is the following six-dimension
system derived at the center of mass of aircraft :\\
\begin{displaymath}\label{equations}
(ED)\quad \quad \left\{\begin{array}{lll}
\dot{V}&=&g\left(\displaystyle\frac{T\cos\alpha-D}{m g}-\sin
\gamma\right)\\
\\
\dot{\gamma}&=&\displaystyle\frac{1}{m V}\left((T\sin\alpha+L)\cos
\mu-m g \cos\gamma \right)\\
\\
\dot{\chi}&=&\displaystyle\frac{(T\sin\alpha+L)\sin \mu}{m V \cos\gamma}\\
\\
\dot{x}&=&V \cos \gamma \cos \chi\\
\\
\dot{y}&=&V \cos \gamma \sin \chi\\
\\
\dot{h}&=&V \sin \gamma\\
\end{array}
\right.
\end{displaymath}
\\
where $T=T(h,V,\delta_x), D=D(h,V,\alpha)$ and $L=L(h,V,\alpha)$.\\
These equations embody the assumptions of a constant weight,
symmetric flight and constant gravitational attraction \cite{Abd07,boif98}.\\
\\
Figure ($1$) shows the forces acting on an aircraft at its center of
gravity during an approach.
\\
\setlength{\unitlength}{0.6cm}
\begin{center}
\begin{picture}(12,7)
\put(0,0){\vector(0,1){7}}
\put(-0.5,6.5){\makebox(0,0){\sffamily\small $h$}}
 \put(0,0){\vector(1,0){11}}
 \put(10.5,-0.5){\makebox(0,0){\sffamily\small $x$}}
\put(6,4){\line(-1,0){4}} \put(6,4){\vector(-1,1){2}}
\put(6,4){\vector(1,1){2}} \put(6,4){\vector(0,-1){2}}
\put(6,4){\vector(-1,-1){2}}
 \put(6,4){\vector(-2,-1){2}}
\put(6,1.6){\makebox(0,0){\sffamily\small W}}
\put(3.9,6.5){\makebox(0,0){\sffamily\small L}}
\put(8.2,6.3){\makebox(0,0){\sffamily\small D}}
\put(3.7,3){\makebox(0,0){\sffamily\small T}}
\put(3.9,1.6){\makebox(0,0){\sffamily\small V}}
\put(3.2,4.3){\makebox(0,0){\sffamily\small HORIZONTAL}}
\put(4.7,3.7){\makebox(0,0){\sffamily\small $\gamma$}}
\put(5.5,3.8){\makebox(0,0){\sffamily\small $($}}
\end{picture}
\begin{figure}[htb]
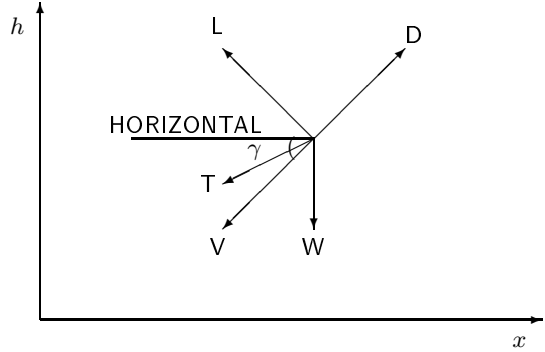

\caption{Aircraft forces during phase of descent}
\end{figure}
\end{center}
Those equations could be applied to conventional aircraft of all
sizes. 
The most dominant aerodynamics affecting results are the lift $L$
and drag $D$, defined as follows \cite{boif98}
\begin{displaymath}
\begin{array}{ll}
L=\frac{1}{2}\rho S V^2 C_{z_{\alpha}}\alpha\\
\\
D=\frac{1}{2}\rho S V^2 \left[C_{x_0}+k_i C_{z_{\alpha}}^2\alpha^2\right] \\
\end{array}
\end{displaymath}
\\
The thrust model chosen, by Matthingly \cite{roux}, depends
explicitly on the aircraft speed, the geometric aircraft height and
the throttle setting
\\
\begin{displaymath}
T=T_0
\delta_x\displaystyle\frac{\rho}{\rho_0}\left(1-M+\displaystyle\frac{M^2}{2}\right)
\end{displaymath}
\\
$T_0$ is full thrust, $\rho_0$ is atmospheric density at the ground
$(=1.225\,\,kg/m^3)$
and $\rho$ is atmospheric density at the height $h \, (\rho=\rho_0(1-22.6\times 10^{-6}h)^{4.26})$.\\
\\
The previous system of equations\, (\ref{equations})\, can be
written in the following generic\\ form :
\begin{displaymath}
\begin{array}{l}
\\
\dot{X}(t)=f(X(t),U(t))\\
\end{array}
\end{displaymath}
where :
\begin{displaymath}
\begin{array}{cccl}
X : &[t_0,t_f]& \longrightarrow & \R^6\\
 & t& \longrightarrow & X(t)=[V(t),\gamma(t),\chi(t),x(t),y(t),h(t)] \;\textrm{are the state variables,}\\
\\
U : &[t_0,t_f]& \longrightarrow & \R^3\\
 & t& \longrightarrow & U(t)=[\alpha(t),\delta_x(t),\mu(t)] \;\textrm{are the control variables,}\\
\end{array}
\end{displaymath}
and $t_0$ and $t_f$ are the initial and final times.\\
\subsection{Constraints}
Search for optimal trajectories minimizing noise must be done in a
realistic flight domain. Indeed, operational procedures are
performed with respect to parameter limits related to the safety of
flight and the
operational modes of the aircraft.\\
\begin{itemize}
\item The throttle stays in some interval
\begin{displaymath}
\delta_{x_{min}}\leq \delta_x \leq \delta_{x_{max}}
\end{displaymath}
\item The speed is bounded
\begin{displaymath}
1.3V_{s_0} \leq V \leq V_{max}
\end{displaymath}
where $V_{s_0}$ the stall velocity, the limited velocity at which
the aircraft can produce enough lift to balance the
aircraft weight. $V_{s_0}$ and $V_{max}$ depend on the type of aircraft. \\
\item The flight path angle providing a measure of the angle of the velocity to the inertial horizontal axis, is bounded
\begin{displaymath}
\gamma_{min}\leq \gamma \leq \gamma_{max}
\end{displaymath}
\item The angle of attack is bounded
\begin{displaymath}
\alpha_{min} \leq \alpha \leq \alpha_{max}
\end{displaymath}
\item The yaw angle and roll angle stay in some prescribed interval
\begin{displaymath}
\begin{array}{lll}
\chi_{min}&\leq \chi& \leq \chi_{max}\\
\mu_{min}&\leq \mu &\leq \mu_{max}\\
\end{array}
\end{displaymath}

\end{itemize}
Those inequality constraints could be formulated as:
\begin{displaymath}
\begin{array}{l}
a \leq C(X(t),U(t)) \leq b
\end{array}
\end{displaymath}
where\,
\begin{equation}\label{contraintes}
\begin{array}{llll}
C :& \R^6 \times \R^3 & \longrightarrow & \R^6\\
&(X(t),U(t))&\longrightarrow&
C(X(t),U(t))=[\gamma(t),V(t),\chi(t),\alpha(t),\delta_x(t),\mu(t)]
\end{array}
\end{equation}
$a$ and $b$ are two constant vectors of $\R^6$.
\subsection{Cost function}
Models and methods used to assess environmental noise problems must
be based on the noise exposure indices used by relevant
international noise control regulations and standards (ICAO
\cite{ICAO93,ICAO95}, Lambert and Vallet \cite{ICAO93}). As
described by Zaporozhest and Tokarev \cite{TZ1}, these indices vary
greatly one from another both in their structure, and in the basic
approaches used in their definitions.\\
\\
The cost function to be minimized may be chosen as any usual
aircraft noise index, which describes the effective noise level of
the aircraft noise event \cite{TZ1,TZ}, like $SEL$ (Sound Exposure
Level), the\,$EPNL$ (Effective Perceived Noise Levels) or the
$L_{eq,\Delta t}$ (Equivalent noise level), ... It is well known
that the magnitude of $L_{eq,\Delta t}$ correlates well with the
effects of noise on human activity, in particular, with the
percentage of highly noise-annoyed people living in regions of
significant aircraft noise impact. This criterion is commonly used,
as basis, for the regulatory basis in many countries. Based on
comparison of noise exposure indices and a comparison of the
methodologies used to calculate the aircraft noise exposure, it can
be concluded that the general form of the most used and accepted
noise exposure index is $L_{eq,\Delta t}$ that we have chosen as an
index (\cite{Jon92,Mont00}). $L_{eq,\Delta T}$ is expressed by:
\begin{equation}\label{Leq}
\begin{array}{ll}
L_{eq,\Delta T}&=10 \log \displaystyle\frac{1}{\Delta
T}\int_{t_0}^{t_f}10^{0.1 L_{P}(t)}dt
\end{array}
\end{equation}
where $t_0$ initial time, $t_f$ final time, $\Delta
T=t_f-t_0$ and $L_{P}(t)$ is the overall sound pressure level (in decibels).\\
\\
We will define in the next subsection the analytic method to compute
the noise level at any reception point.\\
\\
{\bf{Calculation method for aircraft noise levels}}\\
\\
The aircraft noise levels $L_P$ at a receiver is obtained by the
following formula based on works \cite{MGM,KZ}:
\begin{equation}\label{4}
L_{P}=L_{ref}-20 \log_{10} R+\Delta_{atm}+\Delta
_{ground}+\Delta_V+\Delta_D+\Delta_f
\end{equation}
where $L_{ref}$ is the sound level at the source, $20 \log_{10} R $
is a correction due to geometric divergence, $\Delta _{atm}$\,is the
attenuation due to atmospheric absorption of sound. The other terms
$\Delta_{ground},\Delta_V,\Delta_D,\Delta_f$ correspond respectively
to the ground effects, correction for the Doppler,
correction for duration emission and correction for the frequency.\\
\\
In this paper, we have used a semi-empirical model to predict noise
generated by conventional-velocity-profile jets exhausting from
coaxial nozzles predicting the aircraft noise levels represented by
the jet noise \cite{STON} which corresponds to the main predominated
noisy source. It is known that jet noise consists of three principal
components. They are the turbulent mixing noise, the broadband shock
associated noise and the screech tones. At the present time, this
first approximation have been used herein. It seems to be correct in
that step of research because the complexity of the problem. Many
studies have agreed with this model and full-scale experimental data
even at high jet velocities in the region near the jet axis.
Numerical simulation of jet noise generation is not straightforward
undertaking. Norum and Brown \cite{Nur93}, Tam and Auriault
\cite{Tam95,Tam95bis,Tam99} had earlier discussed some of the major
computational difficulties anticipated in such effort. At the
present time, there are reliable to jet noise prediction. However,
there is no known way to predict tone intensity and directivity;
even if it is entirely empirical. This is not surprising for the
tone intensity which is determined by the nonlinearities of the
feedback loop. Obviously, to complete this study we will need to
integrate other noise source models in particular aerodynamics.
\\ \\
Although the numerous aspects of the mechanisms of noise generation
by coaxial jets are not fully understood, the necessity to predict
jet noise has led to the development of empirical procedures and
methods. During the descent phase, the jet aircraft noise as well as
propeller aircraft noise is approximately omni-directional and the
noise emission is decreasing with decreasing speed when assuming
that the power setting is constant. The jet noise results from the
turbulence created by the jet mixing with the surrounding air. Jet
mixing noise caused by subsonic jets is broadband in nature (its
frequency range is without having specific tone component) and is
centered at low frequencies. Subsonic jets have additional shock
structure-related noise components that generally occur at a higher
frequencies. The prediction of jet noise is extremely complex. The
used methods in system analysis and in the engine design usually
employed simpler or semi-empirical prediction techniques. By
replacing the predicted jet noise level \cite{STON} in (\ref{4}), we
obtain the following expression
\begin{equation}\label{niveau}
L_{P}(t)=\left\{\begin{array}{ll}
&141+10\log_{10}\left(\displaystyle\frac{\rho_1}{\rho}\right)^w+
10\log_{10}\left(\displaystyle\frac{V_e}{c}\right)^{7.5}+3\log_{10}\left(\displaystyle\frac{2
s_1}{\pi
d^2}+0.5\right) \\ \\
&+
10\log_{10}\left(\left(1-\displaystyle\frac{v_2}{v_1}\right)^{me}+1.2
\displaystyle\frac{\left(1+\displaystyle\frac{s_2v_2^2}{s_1v_1^2}\right)^4}{\left(1+\displaystyle\frac{s_2}{s_1}\right)^3}\right)+10\log
_{10}s_1\\ \\
&+\log_{10}\displaystyle\frac{\tau_1}{\tau_2}+10\log_{10}\left
(\left(\displaystyle\frac{\rho}{\rho_{ISA
}}\right)^{2}\left(\displaystyle\frac{c}{c_{ISA}}\right)^{4}\right )
\\ \\
&-20\log_{10}R-15\log_{10}(C_D(M_c,\theta))-10\log_{10}(1-M\cos\theta),
\end{array}\right.
\end{equation}
where
\begin{displaymath}
\begin{array}{rll}
v_1&&\textrm{speed of jet gas at inner contours}\\
v_2 &&\textrm{speed of jet gas at outer contours}\\
s_1&&\textrm{area of coaxial engine nozzle at inner contours}\\
s_2 &&\textrm{area of coaxial engine nozzle at outer contours}\\
\tau_1& &\textrm{temperature at inner contours}\\
\tau_2 &&\textrm{temperature at outer contours}\\
\rho_1 &&\textrm{atmospheric density at inner contours}\\
\rho_{ISA}&&\textrm{International Standard Atmosphere density}\\
c_{ISA}&&\textrm{International Standard Atmosphere for speed of sound.}\\
\end{array}
\end{displaymath}
The effective speed $V_e$ is defined by:
$$V_e=v_1[1-(V/v_1)\cos(\alpha_T)]^{2/3}.$$
The angle of attack $\alpha_T$, the upstream axis of the jet
relative to the direction of aircraft motion has been neglected in
this study.\\
The distance source to the observer is
\begin{displaymath}
\begin{array}{l}
R=(x-x_{obs})^2+(y-y_{obs})^2+h^2
\end{array}
\end{displaymath}
where $(x_{obs},y_{obs})$ are the coordinates of the observer.
$\Delta_V$ is expressed by
\begin{displaymath}
\begin{array}{l}
\Delta_V=-15\log(C_D(M_c,\theta))-10\log(1-M\cos\theta)
\end{array}
\end{displaymath}
where $C_D(M_c,\theta)$ indicate the Doppler convection factor:
$$C_D(M_c,\theta)=[(1+M_{c}\cos\theta)^2+0.04 M_c^2],$$
the Mach number of convection is:
$$M_c=0.62(v_1-V)/c$$
and $w
=\displaystyle\frac{3({V_{e}/c})^{3.5}}{0.6+(V_{e}/c)^{3.5}} - 1$.\\
\\
The validity of this improved prediction model is established by
fairly extensive comparisons with model-scale static data
\cite{Tam}. Insufficient appropriate simulated-flight data are
available in the open literature, so verification of flight effects
during aircraft descent has to be established. Analysis by Stone and
al. \cite{STON} has shown that measured data are used to calibrate
the behavior of the used jet noise model and its implications from
theory. Nowadays, the above formulation is being considered
realistic compared to others described in applied and fundamental
literature dealing with jet noise of aircraft during descent
operations.
\\ \\
Taking account the formulas (\ref{Leq}) and (\ref{niveau}), we
obtain our cost function in the following integral function form
\begin{displaymath}\label{cout}
\begin{array}{cll}
J :\mathcal{C}^{1}([t_0,t_f],\R^6)\times
\mathcal{C}^{1}([t_0,t_f],\R^3) &\longrightarrow& \R\\
(X(t),U(t)) &\longrightarrow &
J(X(t),U(t))=\displaystyle\int_{t_0}^{t_f}\ell(X(t),U(t))dt
\end{array}
\end{displaymath}
$J$ is the criterion for optimizing the noise level at the reception
point. It doesn't depend of $\chi,\;\gamma$ and $U$.

\subsection{Optimal Control Problem}

Finding an optimal trajectory, in term of minimizing noise emission
during a descent, can be mathematically stated as an ODE optimal
control problem. We opted different notations $z \equiv X, u \equiv
U$ and $t_0=0$:
\begin{equation}\tag{$OCP$}\label{commande}
\left\{\begin{array}{l}
min \ J(z,u)=\displaystyle\int_{0}^{t_f}\ell(z(t),u(t))dt\\
\\
\dot{z}(t)=f(z(t),u(t)),\forall t \,\in [0,t_f]\\
\\
z_{I_1}(0)=c_1,\,z_{I_2}(t_f)=c_2\\
\\
a \leq C(z(t),u(t))\leq b
\end{array}\right.
\end{equation}

where  $J:\R^{n+m} \rightarrow \R, f :\R^{n+m} \rightarrow \R^n$ and
$ C :\R^{n+m} \rightarrow \R^q$ correspond respectively to the cost
function, the dynamic of the problem, and the constraints defined in
the previous section. The initial and final values for the sate
variables $((x(0),y(0),h(0),V(0))$ and $(x(t_f), y(t_f),h(t_f))$ are
fixed. $n_{\ell} :=|I_1|+|I_2|$
is the total number of fixed limit values of the state variables. \\
\\
To solve our problem (find the optimal control $u(t)$ and the
corresponding optimal state $z(t)$), we discretize the control and
the state with identical grid and transcribe optimal control problem
into nonlinear problem with constraints. The next section may be
helpful in telling how the problem could be solved. They will
present theoretical consideration and computational process that
yield to flight paths minimizing noise levels at a given receiver.
\\
\section{Discrete Optimal Control Problem}

To solve $(OCP)$ different methods and approaches can be used
\cite{BBHLT,W}. In this paper, we use a direct optimal control
technique : we first discretize $(OCP)$ and then solve the resulting
nonlinear programming problem.

\subsection{Discretization}

We use an equidistant discretization of the time interval as
\begin{displaymath}\label{lou}
0=t_0<...<t_{N}=t_f
\end{displaymath}
where :
\begin{displaymath}
t_k=t_0+k {\rm h},\;\;k=0,...,N\;\;\textrm{and}\;\;{\rm
h}=\displaystyle\frac{t_f-t_0}{N}.
\end{displaymath}
\\
Then we consider that $u(.)$ is parameterized as a piecewise
constant function :
\begin{displaymath}
u(t):=u_k\;\;\textrm{for}\;\;t \in [t_{k-1},t_k[
\end{displaymath}
\\
and use a Runge-Kutta scheme (Heun) to discretize the dynamic :
\begin{displaymath}
\left\{\begin{array}{l}
z_{k+1}=z_k+{\rm h}\displaystyle\sum_{i=1}^{s}b_i f(z_{ki},u_{k})\\
z_{ki}=z_k+{\rm h}\displaystyle\sum_{j=1}^{s}a_{ij} f(z_{kj},u_{k}) \\
k=0,\ldots, N-1, \quad i=1,\ldots, s.\\
\end{array}\right.
\end{displaymath}
\\
The new discrete objective function is stated as :
\begin{displaymath}
\tilde J=\displaystyle\sum_{k=0}^{N}\ell(z_k,u_k).
\end{displaymath}
The continuous optimal control problem $(OCP)$ is replaced by the
following discretized control problem :
\begin{displaymath}
(NLP)\quad\quad \left\{
\begin{array}{l}
\min \displaystyle\sum_{k=0}^{N}\ell(z_k,u_k)\\
z_{k+1}=z_k+{\rm h}\displaystyle\sum_{i=1}^{s}b_i f(z_{ki},u_{k}),\quad k=0,\ldots, N-1\\
z_{ki}=z_k+{\rm h}\displaystyle\sum_{j=1}^{s}a_{ij} f(z_{kj},u_{k}), \quad k=0, \ldots, N-1\\
z_{0_{I_1}}=c_1,\;\;z_{N_{I_2}}=c_2\\
\\
a\leq C(z_k,u_k)\leq b,\quad k=0,\ldots, N\\
\end{array}
\right.
\end{displaymath}
\\
\\
To solve (NLP) we developped an AMPL \cite{AMPL} model and used a
robust interior point algorithm KNITRO \cite{Kni}. We choose this
NLP solver after numerous comparisons with some other standard
solvers available on the NEOS (Server for Optimization) platform.


\section{Numerical results}

For different cases and configurations, we consider an aircraft
approach with an initial condition $(x_0=0;y_0=0;h_0=3500\,m)$ a
final condition $(x_f=60000\;m;y_f=5000\;m; h_f=500\;m)$ and for
a fixed $tf=10$ min and a discretization parameter $N=100$ or $N=200$.\\
We first consider the simplest configuration of one single observer
and no additional constraint.\\

\subsection{One fixed observer}

For various positions $(x_{obs},y_{obs})$ of an observer on the
ground (near the aircraft trajectory) we calculate the optimal noise
level $J$ (corresponding to our optimal trajectory $Tr$) and the
noise level $J_1$ corresponding to the trajectory $Tr_1$ that
minimizes the "fuel consumption" (re minimizing) the simple
following model of consumption \cite{boif98}: \\
\begin{displaymath}
CO(h,V,\delta_x)=\int_0^{t_f} C_{SR}T(t)dt
\end{displaymath}
where $C_{SR}$ is supposed constant.\\
\\
\\
The following table summarizes the obtained results.
\begin{table}[htb]
\centering
\begin{displaymath}
\begin{array}{|c|c|c|c|c|c|c|}
\hline
(x_{obs},y_{obs})&J\,(d B)& max(f.e,o.e)&\textrm{CPU}\, (s) & J_1 (d B)  &J_1-J (d B)&\%CO\\
\hline \hline
(0,0)& 45.92&7.92e-07&10&47.03&1.1 & 36\%\\
\hline (0,2500) & 44.95&8.76e-07&8.4 &46.32&1.4&36\%
 \\
\hline (0,5000) & 43.27&3.44e-07&9.8&44.93&1.7&36\%
\\
\hline\hline
(20000,0)& 48.97&4.66e-07&10.4& 51.04&2.10&33\%\\
\hline
(20000,2500) & 49.58&2.81e-08 &10.6& 51.55&2& 26\%\\
\hline
(20000,5000) & 47.18& 8.36e-07&9.8&49.70&2.5&36\%\\
\hline\hline (40000,0)&47.59&6.73e-07&14.6&49.52 &2&26\%
\\
\hline
(40000,2500) & 50.76&4.21e-07&6.9&52.87 &2.11&26\%\\
\hline (40000,5000) & 49.74&6.92e-07&11.9&51.72& 2&24\%
\\
\hline\hline (60000,0)&42.85&5.29e-07&7.2&45.00&2.15&34\%
\\
\hline (60000,2500) & 45.04&6.90e-07&6.7&48.014&3&36\%
\\
\hline (60000,5000) & 48.84&7.10e-07&6&54.18&5.34& 38\%
 \\
\hline
\end{array}
\end{displaymath}
\caption{Noise minimization}\label{kk}
\end{table}
\\
For each case, the algorithm (KNITRO\cite{Kni}) found a solution
with a very high accuracy. The computation of $Tr_1$ have been done
only one time; it needs $7 s$ with an accuracy of
$max(f.e,o.e)=4.24e-07$.\\
\\
The third column of Table $1$ measure the maximum of feasibility
error and optimality error, the fourth one gives an idea on the
computation effort (namely the CPU time). The two last columns
correspond to the noise
reduction and $\%$ of exceeded consumption : $\%CO=(CO(Tr)-CO(Tr_1))/CO(Tr)$.\\
Our trajectory that minimizes the noise consume about $31\%$ more
than the trajectory minimizing the consumption. \\
\\
\\
The following figure showing the solution trajectory $Tr$, where the
fixed observer presents a certain area near the airport :
\newpage
\begin{figure}[htb]\label{F1}
\begin{center}
\includegraphics[scale=0.7]{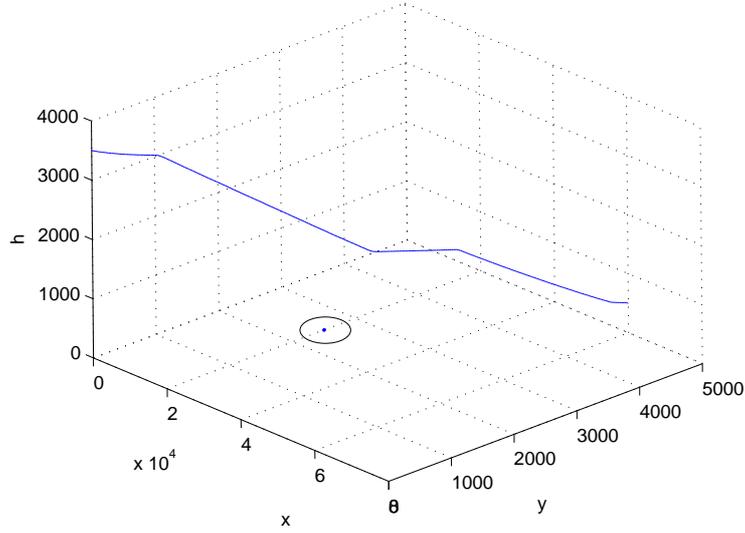}
\caption{Trajectory in 3D}
\end{center}
\end{figure}

The following figures present the state and control variables of the
optimal trajectory $Tr$.
\begin{figure}[htb]\label{F2}
\begin{center}
\includegraphics[scale=0.7]{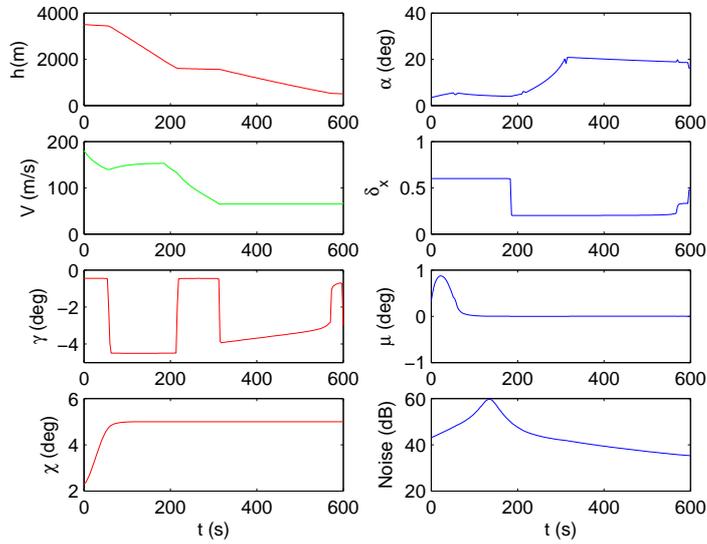}
\caption{Solutions of $(NLP)$}
\end{center}
\end{figure}
We remark that the optimal variables $h,V,\chi,\alpha$ and $\mu$
present some large constant stages, while $\gamma$ and $\delta_x$
are bang-bang.

\subsubsection{One fixed observer with an additional consumption
constraint}

Table (1) shows that the optimal trajectory $Tr$ consumes about
$31\%$ more that $Tr_1$. This fact makes of interest some additional
constraint on the consumption.
\newpage
We define a new problem :
\begin{equation} (OCP)_2\quad \label{commande}
\left\{\begin{array}{l}
\min\, L_{eq,\Delta t}\\
\dot{z}(t)=f(z(t),u(t))\\
z_{I_1}(0)=c_1,\,z_{I_2}(t_f)=c_2\\
CO(z(t),u(t))<= 1.1 CO(Tr_1)\\
\end{array}\right.
\end{equation}
\\
This problem can be solved using the same techniques
(discretization,...) and the same configurations. We obtain the
following results.\\
\begin{displaymath}
\begin{array}{|c|c|c|c|c|c|c|c|c|}
\hline
(x_{obs},y_{obs})&J\,(d B)& max(f.e,o.e) &J_1 (d B)&\textrm{CPU}\, (s) \\
\hline \hline (0,0)&46.32&9.94e-07&47.03&9.6 \\

\hline (0,2500) &45.44&9.17e-07&46.32&13.2
\\ \hline (0,5000) &43.94&9.49e-07&44.93&13.1
\\
\hline\hline (20000,0)&49.79&7.25e-07&51.04&18.6
\\
\hline (20000,2500)& 50.32&8.75e-07&51.55&8.5
\\
\hline (20000,5000) & 48.33&9.37e-07&49.70&20
\\
\hline\hline (40000,0)&48.16&7.54e-07&49.52&33.3
\\
\hline (40000,2500) & 51.43&6.55e-07&52.87&9.3
\\ \hline (40000,5000) &50.31&6.90e-07&51.72&20.4
\\
\hline\hline (60000,0)&43.81& 9.74e-07&45.00&34.8
\\
\hline (60000,2500) & 46.40&7.69e-07&48.014&27.9
\\
\hline (60000,5000) &51.32&9.30e-07&54.18&17.8
 \\
\hline
\end{array}
\end{displaymath}
\begin{table}[htb]
\centering
\caption{Noise minimization}\label{kk}
\end{table}
\\
\\
The following figures present the state and control variables and
the solution trajectory $Tr$ :
\begin{figure}[htb]\label{shema3}
\begin{center}
\includegraphics[scale=0.7]{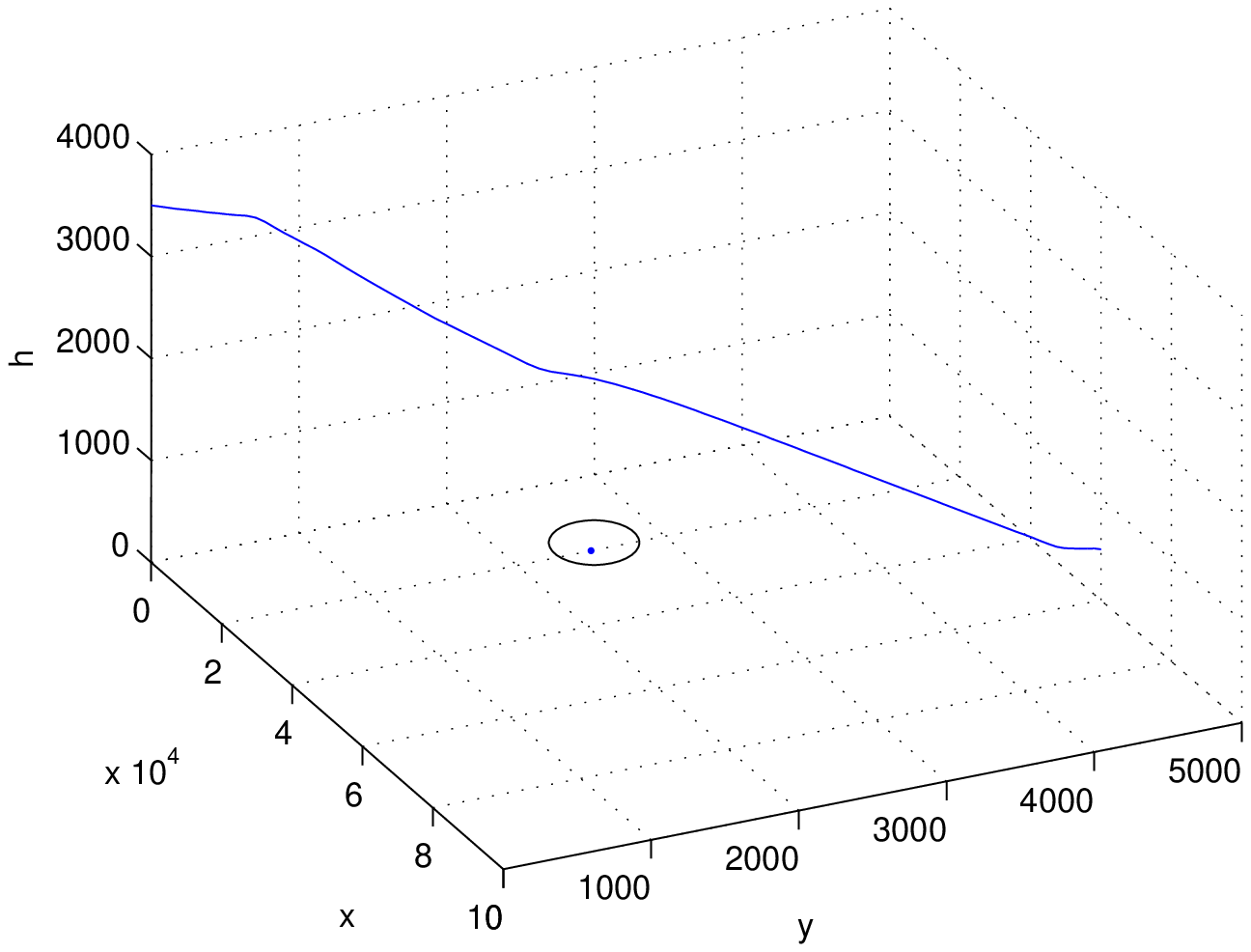}
\caption{Trajectory in 3D}
\end{center}
\end{figure}
\newpage
\begin{figure}[htb]\label{shema4}
\begin{center}
\includegraphics[scale=0.8]{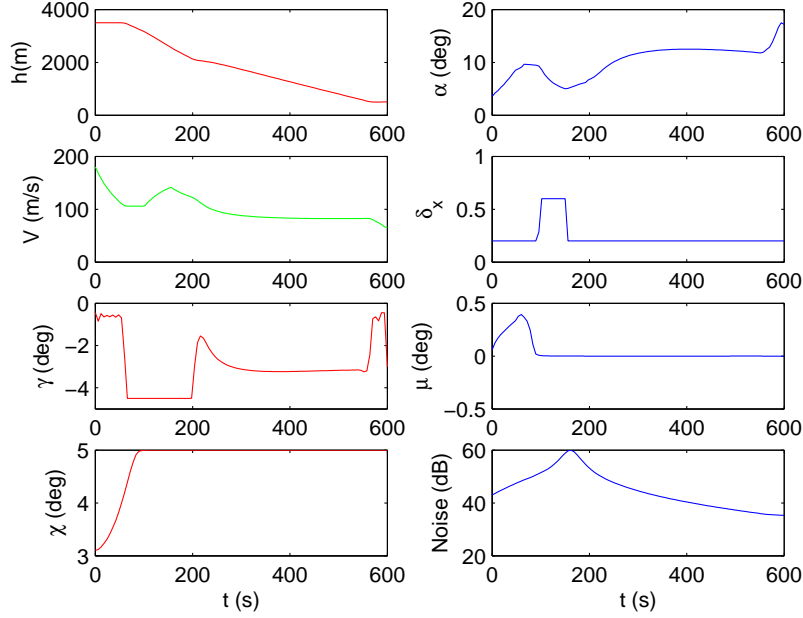}
\caption{Solutions of $(OCP)_2$}
\end{center}
\end{figure}

We obtain approximately the same characteristic for the trajectory
while the noise reduction is still significative.\\
\\
\\
\subsection{Several observers fixed on the ground}
We can generalize the processus of minimization for several
observers. In this case, we
minimize the maximum of noise corresponding to several observers.\\
The problem to solve is written as follows :\\
\begin{equation}\quad (OCP)_3\label{commande}
\left\{\begin{array}{l}
\min \vartheta \\
\vartheta >= J_{{obs}_j}\\
\dot{z}(t)=f(z(t),u(t))\\
z_{I_1}(0)=c_1,\,z_{I_2}(t_f)=c_2\\
\end{array}\right.
\end{equation}
\\
where $J_{{obs}_j}$ are the noise levels corresponding to $j$
fixed observers.\\

Once again, we use a direct method to solve the problem $(OCP)_3$ with the same modeling language and software. \\
We choose five observers : $(obs_1(0,0)$, $obs_2(20000,2500)$,
$obs_3(40000,5000)$, \\$obs_4(60000,0)$, $obs_5(60000,5000))$ which
represent a certain area near the airport. We obtained an optimal
solution, the obtained noise level is $49.65$ dB, the accuracy of
the results is still very high $(max(f.o,e.o)=5.89e-07)$ and the
algorithm takes no more $130 s$ on a standard PC. This trajectory is about 5 dB less that $J_1 (54.18$ dB).\\
The trajectory characteristics are given in the following figures:
\newpage
\begin{figure}[htb]
\begin{center}
\includegraphics[scale=0.75]{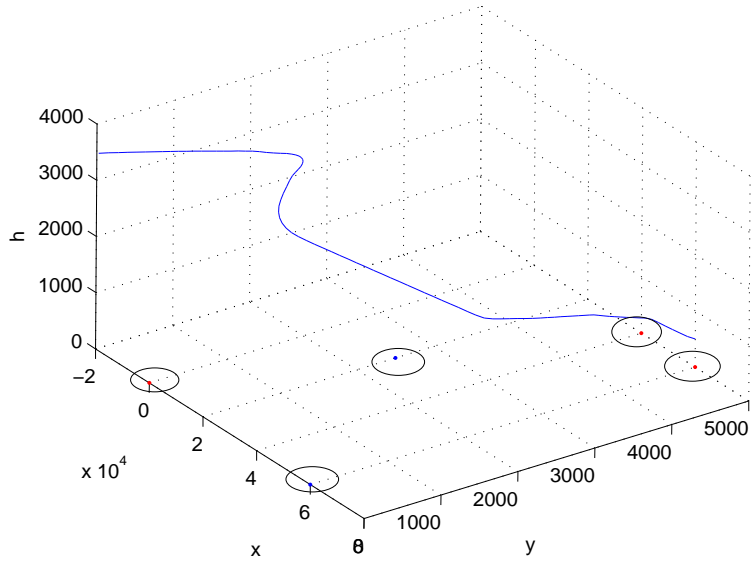}
\caption{Trajectory in 3D with several observers}
\end{center}
\end{figure}
\begin{figure}[htb]\label{shema}
\begin{center}
\includegraphics[scale=0.75]{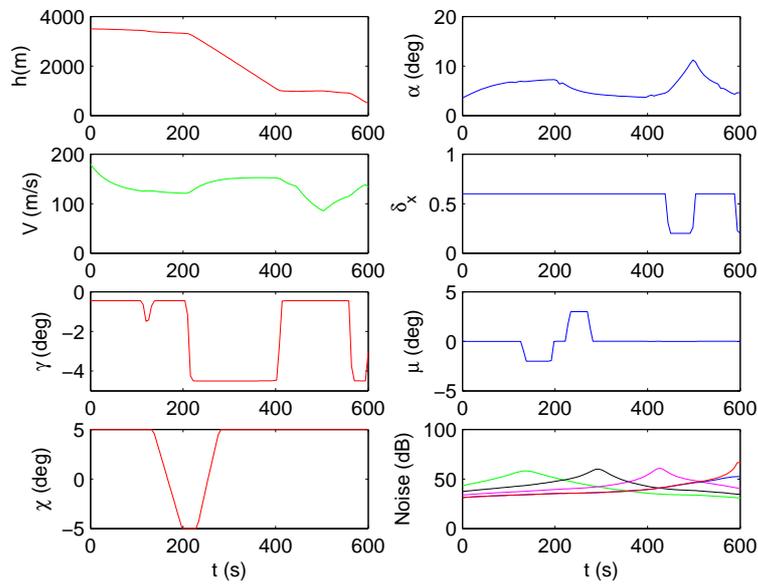}
\caption{Solutions of $(OCP)_3$}
\end{center}
\end{figure}
Almost all state and control variables (except $\alpha$ and $V$)
present large constant stages. The control variables $\delta_x$ and
$\mu$ are bang-bang between their prescribed bounds.\\
\newpage
\section{Conclusion}

We have performed a numerical computation of the optimal control
issue. An optimal solution of the discretized problem is found with
a very high accuracy. A noise reduction is obtained during the phase
approach by considering the configuration of one and several
observers.
The trajectory obtained presents interesting characteristics and performances.\\
\\
Extensions of the analysis on the current problem should include
other source of noise. This feature is particularly important since
improved noise model that better represent individual noise sources
(engines, airframe,...). It should be remembered that this model is
focused on single event flight. Additional researches are needed to
fully assess the influence of wind and other atmospheric conditions
on noise prediction process.
The noise studies have, as yet, been limited to a single aircraft
type equipped with two engines. Further research should consider
multiple flights model considering airport
capacity and nearby configuration.\\

{\bf{Acknowledgments}}\\
\\
We would like to thank Thomas Haberkorn for his constructive
comments and his valuable help  and suggestions for computational
aspects.

\end{document}